\documentclass[english]{article}
\usepackage[latin9]{inputenc}
\usepackage{textcomp}
\usepackage{amstext}
\usepackage{setspace}

\makeatletter

\usepackage{amstext}

\newcommand{\lyxmathsym}[1]{\ifmmode\begingroup\def\b@ld{bold}
  \text{\ifx\math@version\b@ld\bfseries\fi#1}\endgroup\else#1\fi}

\usepackage{babel}

\makeatother

\usepackage{babel}
\begin{document}

\title{On Stanley's theorem \& It's Generalizations}

\author{Suprokash Hazra}
\maketitle
\begin{abstract}
\begin{doublespace}
In this paper we go on to discuss about Stanley's theorem in Integer
partitions. We give two different versions for the proof of the generalization
of Stanley's theorem illustrating different techniques that may be
applied to profitably understand the underlying structure behind the
theorem.\end{doublespace}

\end{abstract}
\begin{doublespace}
 \vspace{5pt}
 
\end{doublespace}

\begin{doublespace}

\subsection*{Introduction:}
\end{doublespace}

\begin{doublespace}
Stanley's theorem is an important result in the theory of partitions.
Various generalizations have been made to it over the course of time
including Elder's theorem{[}1{]} and Dastidar \& Gupta's{[}2{]} subsequent
generalization.In their paper they go on to show how the sum of the
number of distinct members of the partition is not just equal to the
number of 1's present in the partition of the same number but is also
equal to the sum of the number of i's in the partitions of all the
numbers from n to (n+i-1), where i can be any positive integer. In
this paper we go on to provide different combinatorial arguments to
prove such generalizations. A proof of the generalization involves
the use of tilings of a 1x$\infty$ board. This concept can be further
generalised to cover similar properties for overpartitions as well.
To prove the generalization we first make use of some lemmas which
follow subsequently in the article.

\vspace{8pt}

\end{doublespace}

\begin{doublespace}

\subsection*{LEMMA 1.1 :}
\end{doublespace}

\begin{doublespace}
Let n and k be two positive integers with $k\leq n$, then for each
positive 
\end{doublespace}

integer i, $(1\leq i\leq n)$, we have $n+i=q_{i}k+r_{i}$, where
$0\le r_{i}<k$, 

\begin{doublespace}
then $q_{i}=q_{0}$ , for all $1\le i\le s-1$ 

$=q_{0}+1$, for all i $\ensuremath{\ge}$ s 

And $r_{i}=r_{0}+i$ , for all 1 $\ensuremath{\le}$i $\ensuremath{\le}$s-1,

$=i-s$ , for all i $\ensuremath{\ge}$s where $s=k-r_{0}$

\vspace{8pt}

\end{doublespace}

\begin{doublespace}

\subsubsection*{Proof of lemma 1.1 : }
\end{doublespace}

\begin{doublespace}
Let us consider the two integers n+i and k,

then by division algorithm there exists two integers $q_{i}$ and
$r_{i}$

such that,$n+i=q_{i}k+r_{i}$ with $0\le r_{i}<k$. Now, $n=q_{0}k+r_{0}$
.

We assume, $k-r_{0}=s$

Then we have the following,

$n=q_{0}k+r_{0};$

$n+1=q_{0}k+(r_{0}+1);$

$n+2=q_{0}k+(r_{0}+2)$

...

$n+s-1=q_{0}k+(r_{0}+s-1)$ 

$n+s=(q_{0}+1)k$ 

$n+s+1=(q_{0}+1)k+1$ 

$n+s+2=(q_{0}+1)k+2$

.....................

$n+k-1=(q_{0}+1)k+(r_{0}-1)$

\vspace{6pt}

By the above formula it is clear that

$q_{i}=q_{0}$ for all $1\leq i\leq s-1$

=$q_{0}+1$ for all $i\geq s$

Hence the lemma 1.1

Notation: $Q_{k}(n+i):=$Number of occurences of k 
\end{doublespace}

in the unordered partitions of $(n+i)$

\begin{doublespace}
$B(n):=$Sum of the total number of distinct members 

in all the partitions of n.

$A(n)=\sum_{i=0}^{k-1}Q_{k}(n+i)$=Sum of all the number of occurrences
of k in all 

unordered partitions of n+i for all $0\ensuremath{\le}i\ensuremath{\le}k-1$
\end{doublespace}

\begin{doublespace}

\subsection*{LEMMA 1.2 : }
\end{doublespace}

\begin{doublespace}
$Q_{k}(n+i)=\sum_{j=1}^{q_{i}}P(n+i-jk)$, for some fixed i.

\vspace{8pt}

\end{doublespace}

\begin{doublespace}

\subsubsection*{Proof of lemma 1.2 : }
\end{doublespace}

\begin{doublespace}
For a fixed i, considering the two integers n+i and k, then by division 

algorithm there exists two integers $q_{i}$ and $r_{i}$ such that,
$n+i=q_{i}k+r_{i}$

with $0\le r_{i}<k$ (as mentioned in lemma 1.1)

For a fixed i, the number of partitions of n+i,

where k occurs exactly once is $P(n+i-k)-P(n+i-2k)$,

the number of partitions of n+i,

where k occurs exactly twice is $P(n+i-2k)-P(n+i-3k)$,

The number of partitions of n+i, where k occurs

exactly $(q_{i}-1)$ times is $P(n+i-(q_{i}-1)k)-P(n+i-q_{i}k)$

And , the number of partitions of n+i, where k occurs 

exactly $q_{i}$ times is $P(r_{i})$.

So $Q_{k}(n+i)$

= The number of occurrences of k in all unordered partitions of n+i

= $P(n+i-k)-P(n+i-2k)$ + $2[P(n+i-2k)-P(n+i-3k)]$ +\dots . 

$...+(q_{i}-1)[P(n+i-(q_{i}-1)k)-P(n+i-q_{i}k)]+q_{i}P(r_{i})$

\vspace{5pt}

=$P(n+i-k)+P(n+i-2k)+P(n+i-3k)+\ldots..$ 

$...+P(n+i-(q_{i}-1)k)-(q_{i}-1)P(n+i-q_{i}k)+q_{i}P(r_{i})$

\vspace{5pt}

=$P(n+i-k)+P(n+i-2k)+P(n+i-3k)+\ldots..$ 

$...+P(n+i-(q_{i}-1)k)-q_{i}P(n+i-q_{i}k)+P(n+i-q_{i}k)+q_{i}P(r_{i})$

\vspace{5pt}

=$P(n+i-k)+P(n+i-2k)+P(n+i-3k)+\ldots..$ 

$.....+P(n+i-(q_{i}-1)k))+P(n+i-q_{i}k)[since,n+i-q_{i}k=r_{i}]$

\vspace{5pt}

=$\sum_{j=1}^{q_{i}}P(n+i-jk)$,

\vspace{5pt}

Hence the lemma 1.2.

\vspace{8pt}

\emph{\large{}Observation:}\emph{ }When $t\varepsilon\{1,2,...,n\}$
is fixed for all $1\ensuremath{\le}$ i $\ensuremath{\le}$n ,we have, 

$\sum_{i=0}^{k-1}[P(n+i-tk)]=P(n-tk)+P(n+1-tk)+P(n+2-tk)+\ldots..$

$+P(n+k-1-tk)$

$=P(n-tk)+P(n-(tk-1))+P(n-(tk-2))+\ldots.+P(n-((t-1)k+1))$.

\vspace{8pt}

\end{doublespace}

Now we come to our main theorem which states that, 

\begin{doublespace}

\subsection*{THEOREM 1 :}
\end{doublespace}

\begin{doublespace}
$A(n)$ = $B(n)$ for all integers $n$
\end{doublespace}

\begin{doublespace}

\subsection*{Proof of theorem: }
\end{doublespace}

\begin{doublespace}
We consider $A(n)=\text{\ensuremath{\sum}}_{i=0}^{k-1}Q_{k}(n+i)$
=$\sum_{i=0}^{k-1}[\sum P(n+i-jk)]$ 
\end{doublespace}

\{By lemma 1.2\}

\begin{doublespace}
=$\sum_{i=0}^{k-1}[P(n+i-k)]$+ $\text{\ensuremath{\sum}}_{i=0}^{k-1}[P(n+i-2k)$
+\dots .

\vspace{4pt}

+ $\sum_{i=0}^{k-1}[P(n+i-(q_{i}-1)k))]+\sum_{i=0}^{k-1}[P(n+i-q_{i}k)]$ 

\vspace{9pt}

from which it follows that $A(n)$

= $[P(n-1)+P(n-2)+P(n-3)+\ldots.+P(n-k)]$ + $[P(n-(k+1))+$ 

\vspace{4pt}

$P(n-(k+2))+\text{\ensuremath{\dots}}.+P(n-2k)]$\dots{} + $[P(n-((q_{0}-1)k+1)+$ 

\vspace{5pt}

$P(n-((q_{0}-1)k+2)+\text{\ensuremath{\dots}}.+P(n-q_{0}k)]$$+\sum_{i=s}^{k-1}[P(n+i-q_{i}k)]$ 
\end{doublespace}

(by above observation)

\begin{doublespace}
\vspace{9pt}

= $P(n-1)+P(n-2)+P(n-3)+\ldots.+P(n-q_{0}k)+\sum_{i=s}^{k-1}P(r_{i})$

(Using lemma 1.1)

\vspace{9pt}

$=P(n-1)+P(n-2)+P(n-3)+\ldots.+P(r_{0})$+ 

\vspace{4pt}

$[P(r_{0}-1)+P(r_{0}-2)+\ldots.+P(1)+1]$

This is the expression for A(n). Now to find the expression for B(n)
using $P(1),P(2),\ldots.,P(n-1)$. For finding the expression for
B(n),we will look at this sum to count the number of partitions of
n in which the number i appears and sum those result for all $1\ensuremath{\le}i\ensuremath{\le}n.$Now
the number of partitions of n in which i appears is $P(n-i)$ for
all $\ensuremath{1\le i\leq n-1}$ and with the special case $i=n$
is 1. Hence $B(n)=1+P(1)+P(2)+\ldots.+P(n-3)+P(n-2)+P(n-1).$
\end{doublespace}

\begin{verse}
Hence A(n)= B(n)

\begin{doublespace}
This completes our proof. \end{doublespace}

\end{verse}
\begin{doublespace}

\section*{2. A TILLING PROOF OF THE EXTENSION OF STANLEY'S THEOREM }
\end{doublespace}

\begin{doublespace}
We consider a 1x$\infty$ board. We will tile this board using white
squares and finitely many black squares. We allow the stacking of
the black squares. Let T be the set of all such tilings. Let a white
tile will have measure 1 and a black tile in position i will have
measure $q^{i}$. Now we define the measure of a tiling to be M(t)=
$\prod$ m(t). (where the product runs over all the squares for a
particular tiling t$\varepsilon$ T. And m denotes the measure of
the squares for a particular tiling t.
\end{doublespace}

\begin{doublespace}

\subsection*{Partition representation of a tiling:}
\end{doublespace}

\begin{doublespace}
Suppose we assume a partition of $m=\lyxmathsym{\textmu}_{1}+\lyxmathsym{\textmu}_{2}+\ldots..\lyxmathsym{\textmu}_{K}$,
denoting this partition by \textmu . Then we associate a tiling $t_{\lyxmathsym{\textmu}}$
such that, it has a black tiles in position $\text{\textmu}_{1},\text{\textmu}_{2}...,\lyxmathsym{\textmu}_{K}$.
Now as the numbers $\text{\textmu}_{1},\text{\textmu}_{2},..,\lyxmathsym{\textmu}_{K}$
may not all distinct so we can have more than one black tiles in some
position. Also we obtain $M(t_{\lyxmathsym{\textmu}})=q^{\text{\textmu}_{1}+\text{\textmu}_{2}+\text{\ensuremath{\dots}}..\text{\textmu}_{K}}$

Now we need these following lemmas in order to prove the theorem. 
\end{doublespace}

\begin{doublespace}

\subsection*{LEMMA 2.1 : }
\end{doublespace}

\begin{doublespace}
For each $1\ensuremath{\le}i\ensuremath{\le}k-1$, the number of partitions
of n with atleast (k-i) times
\end{doublespace}

r = Number of partitions of n+ir with atleast r times k. 

\begin{doublespace}

\subsubsection*{Proof of lemma 2.1 :}
\end{doublespace}

\begin{doublespace}
Consider A= Set of all tilings with atleast (k-i) black tiles in position
r.

And B=Set of all tilings with atleast r black tiles in position k.

Now we define T: A\textrightarrow B, by the following way,

We take a tiling from A, then we remove (k-i) black tiles from position
r and 
\end{doublespace}

add r black tiles in position k.

\begin{doublespace}
Then T is well defined map. In order to prove T is bijective

we define a mapping S:B\textrightarrow A, by the following way,we
take a tiling from B,

then we remove r black tiles from position k and add (k-i) black tiles
in

position r. Then S is well defined map. Also S T = Id$_{A}$ and T
S = Id$_{B}$,

the identity mappings on A and B respectively.

Now clearly the change of measure under the map T is $q^{-r(k-i)+kr}=q^{ir}$.

Now take A{*}=Set of all tilings of A with measure $q^{n}$.

Then for each tiling of A{*} we obtain a bijective Correspondence
with a tiling
\end{doublespace}

of B of measure $q^{n+ir}$.

\begin{doublespace}
Now every tiling gives a partitions and vice versa.

So by the Partition representation of a tiling,

we have that, the number of partitions of n

with atleast (k-i) times r = number of partitions of n+ir 
\end{doublespace}

with atleast r times k.

\begin{doublespace}
Hence the lemma 2.1. 
\end{doublespace}

\begin{doublespace}

\subsection*{LEMMA 2.2 : }
\end{doublespace}

\begin{doublespace}
Number of partitions of n with atleast i times r = The number of partitions
of n with atleast r times i.
\end{doublespace}

\begin{doublespace}

\subsubsection*{Proof of lemma 2.2:}
\end{doublespace}

\begin{doublespace}
Consider M= Set of all tilings with atleast i black tiles in position
r. And N=Set of all tilings with atleast r black tiles in position
i. Now we define Q: M$\rightarrow$N, by the following way, We take
a tiling from M, then we remove i black tiles from position r and
add r black tiles in position i. Then Q is well defined map. In order
to prove Q is bijective we define a mapping Z:N$\rightarrow$M, by
the following way,we take a tiling from N, then we remove r black
tiles from position i and add i black tiles in position r. Then Z
is well defined map. Also QZ= Id$_{N}$ and ZQ = Id$_{M}$, the identity
mappings on N and M respectively. Now clearly there is no change of
measure under the map Q. Now take M{*}=Set of all tilings of A with
measure $q^{n}$. Then for each tiling of A{*} we obtain a bijective
Correspondence with a tiling of N of measure $q^{n}$ . Now every
tiling gives a partitions and vice versa. So by the Partition representation
of a tiling, we have the lemma 2.2. 
\end{doublespace}

\begin{doublespace}

\subsection*{LEMMA 2.3 : }
\end{doublespace}

\begin{doublespace}
The number of partitions of n with atleast 1 times r = The number
of partitions of (n+kj-r) with atleast j times k. we assume $n=qk+r$,and
$1\ensuremath{\le}j\ensuremath{\le}q$. In particular,$V_{j}^{k}(n+kj-r)=V_{s}^{k}(n+ks-r)$,
for all $1\le j,s\le n$ Where $V_{p}^{q}(m)$= The number of partitions
of m with atleast p times q. 

\vspace{8pt}

\end{doublespace}

\begin{doublespace}

\subsubsection*{Proof of lemma 2.3:}
\end{doublespace}

\begin{doublespace}
Let C = Set of all tilings with atleast 1 black tiles in position
r. And D= Set of all tilings with atleast j black tiles in position
k. Now we define F: C\textrightarrow D, by the following way, We take
a tiling from C, then we remove 1 black tiles from position r and
add j black tiles in position k. Then F is well defined map. In order
to prove F is bijective we define a mapping G:D\textrightarrow C,
by the following way, we take a tiling from D, then we remove j black
tiles from position k and add 1 black tiles in position r. Then G
is well defined map. Also GF = Id$_{C}$and FG = Id$_{D}$, the identity
mappings on C and D respectively. Now clearly the change of measure
under the map F is $q^{kj-r}$.

$C^{\#}$=Set of all tilings of C with measure $q^{n}$. Then for
each tiling of $C^{\#}$ we obtain a bijective Correspondence with
a tiling of D of measure $q^{n+kj-r}$.

Now every tiling gives a partition and vice versa. So by the Partition
representation of a tiling, we have that, the number of partitions
of n with atleast 1 times r = number of partitions of n+kj-r with
atleast j times k.

So notationally, $V_{j}^{k}(n+kj-r)=V_{1}^{r}(n)$.

Since the right hand side of the above expression does not depend
on j,so clearly we have that, $V_{j}^{k}(n+kj-r)=V_{s}^{k}(n+ks-r)$,
for all $1\ensuremath{\le}j,s\ensuremath{\le}n$

Hence the lemma 2.3.
\end{doublespace}

\begin{doublespace}

\end{doublespace}

\begin{doublespace}
Now we come to our main theorem.
\end{doublespace}

\begin{doublespace}

\subsection*{Proof of theorem :}
\end{doublespace}

\begin{doublespace}
First we put r=1 in lemma 2.1.

Then we get, the number of partitions of n with atleast (k-i) times
1 = number of partitions of n+i with atleast 1 times k.

\vspace{5pt}

Notationally $V_{(k-i)}^{1}(n)=V_{1}^{k}(n+i)$ for all 1$\ensuremath{\le}i\ensuremath{\le}k-1$
We sum this over all $1\ensuremath{\le}i\ensuremath{\le}k-1$ and
obtain, $\sum_{i=1}^{k-1}V_{(k-i)}^{1}(n)=\sum_{i=1}^{k-1}V_{1}^{k}(n+i)$
\dots \dots{} (1)

\vspace{5pt}

Now we put r=1 in lemma 2.2. and obtain, $V_{i}^{1}(n)=V_{1}^{i}(n)$
We sum this over all k $\ensuremath{\le}i\ensuremath{\le}n$ and obtain,
$\sum_{i=k}^{n}V_{i}^{1}(n)=\sum_{i=k}^{n}V_{1}^{i}(n)$ \dots ..(2)

\vspace{5pt}

We add (1) and (2),we obtain

$\sum_{i=1}^{k-1}V_{(k-i)}^{1}(n)+\sum_{i=k}^{n}V_{i}^{1}(n)=\sum_{i=1}^{k-1}V_{1}^{k}(n+i)+\sum_{i=k}^{n}V_{1}^{i}(n)$

\vspace{5pt}

$\sum_{i=1}^{n}V_{i}^{1}(n)=\sum_{i=1}^{k-1}V_{1}^{k}(n+i)+\sum_{i=k}^{n}V_{1}^{i}(n)$

Number of 1's present in the all unordered 

partitions of n
\end{doublespace}

= $\sum_{i=1}^{k-1}V_{1}^{k}(n+i)+\sum_{i=k}^{n}V_{1}^{i}(n)$

\begin{doublespace}
\vspace{5pt}

= $\sum_{i=1}^{k-1}V_{1}^{k}(n+i)+V_{1}^{k}(n)+\sum_{i=k+1}^{n}V_{1}^{i}(n)$

\vspace{5pt}

= $\sum_{i=0}^{k-1}V_{1}^{k}(n+i)+\sum_{i=k+1}^{n}V_{1}^{i}(n)$\dots \dots \dots \dots \dots{}
(4)

\vspace{5pt}

\end{doublespace}

Now by lemma 2.3 

\begin{doublespace}
\vspace{5pt}

\end{doublespace}

$\sum_{i=k+1}^{n}V_{1}^{r}(n)=\sum_{i=k+1}^{n}V_{j}^{k}(n+kj-r)$

\begin{doublespace}
\vspace{5pt}

$\sum_{i=k+1}^{n}V_{1}^{r}(n)=\sum_{i=k+1}^{2k}V_{j}^{k}(n+kj-r)+\sum_{i=2k+1}^{3k}V_{j}^{k}(n+kj-r)+\ldots\ldots$
$+\text{\ensuremath{\sum}}_{i=(q-1)k+1}^{qk}(V_{j}^{k}(n+kj-r))$

\vspace{5pt}

$\sum_{i=k+1}^{n}V_{1}^{r}(n)=\sum_{i=k+1}^{2k}V_{2}^{k}(n+2k-r)+\sum_{i=2k+1}^{3k}V_{3}^{k}(n+3k-r)+\ldots\ldots$
$\text{\ensuremath{\dots}}.\text{\ensuremath{\sum}}_{i=(q-1)k+1}^{qk}(V_{q}^{k}(n+kq-r))$

\vspace{5pt}

$\sum_{i=k+1}^{n}V_{1}^{r}(n)=\sum_{i=0}^{k-1}V_{2}^{k}(n+i)+\sum_{i=0}^{k-1}V_{3}^{k}(n+i)+\ldots\ldots\ldots\ldots$\dots{} 

$..+\text{\ensuremath{\sum}}_{i=0}^{k-1}V_{q}^{k}(n+i)$\dots \dots \dots{}
(6)

\vspace{8pt}

Now by (4) and (6) we have that, Number of 1's present in the all
unordered partitions of n

\vspace{5pt}

$=\sum_{i=0}^{k-1}V_{1}^{k}(n+i)+\sum_{i=0}^{k-1}V_{2}^{k}(n+i)+\sum_{i=0}^{k-1}V_{3}^{k}(n+i)+\ldots\ldots\ldots\ldots$ 

$..+\text{\ensuremath{\sum}}_{i=0}^{k-1}V_{q}^{k}(n+i)$

\vspace{5pt}

$=\sum_{j=1}^{q}V_{j}^{k}(n)+\sum_{j=1}^{q}V_{j}^{k}(n+1)+\ldots\ldots\ldots.+\sum_{j=1}^{q}V_{j}^{k}(n+k-1)$ 

\vspace{5pt}

$=Q_{k}(n)+Q_{k}(n+1)+Q_{k}(n+2)+\text{\ensuremath{\dots\dots\dots\dots}}.+Q_{k}(n+k-1)$

\vspace{5pt}

$=\sum_{i=0}^{k-1}Q_{k}(n+i)$ 

\vspace{5pt}

This completes our proof. 
\end{doublespace}

\begin{doublespace}

\section*{Conclusion:}
\end{doublespace}

\begin{doublespace}
In this paper we have tried to demonstrate how simple combinatorial
proofs can be effectively applied to problems involving unordered
partitions. Similar to the extension of Stanley\textquoteright s theorem,
we can apply similar techniques to prove further extensions of Elder\textquoteright s
theorem and not just for regular partitions but for all classes of
overpartitions. The concept of tilings could also pe profitably extended
to study similar results in the area of planar partition 
\end{doublespace}


\begin{thebibliography}{1}
\bibitem{key-1}Wolfram Mathworld. Elder's Theorem. 2010. Available
at http://mathworld. wolfram.com/EldersTheorem.html

\bibitem{key-2}Dastidar, Gupta. Generalization of a few results in
integer partitions

\bibitem{key-3}Andrews, George, The Theory of Partitions, in ``Encyclopedia
of Mathematics and Its Applications,\textquotedblright{} Vol. 2, Addison-Wesley,
Reading, MA, 1976.\end{thebibliography}
\end{document}